\documentclass[12pt,leqno,twoside,a4paper]{article}
\usepackage{amsfonts,amsmath,amsthm,amssymb}
\usepackage[english]{babel}
\newtheorem{theorem}{Theorem}[section]
\newtheorem{prop}[theorem]{Proposition}
\newtheorem{lemma}[theorem]{Lemma}
\newtheorem{corollary}[theorem]{Corollary}
\theoremstyle{definition}
\newtheorem{definition}[theorem]{Definition}
\newtheorem{example}[theorem]{Example}
\theoremstyle{remark}
\newtheorem{remark}[theorem]{Remark}
\newtheorem{Problem}{\bf Problem}
 \numberwithin{equation}{section}

\newcommand{\si}{\sigma}
\newcommand{\p}{\phi_\omega}
\newcommand{\om}{\omega}
\newcommand{\vm}{\overline{\omega}}
\newcommand{\de}{\delta}
\newcommand{\ov}[1]{\overline{#1}}
\newcommand{\dw}{\delta_{\omega}}
\newcommand{\T}{R[x;\sigma ,\delta_{\omega}]}
\newcommand{\Tn}{R[x;\sigma ,\delta]}
 \newcommand{\an}[2]{\mbox{\rm r.ann}_{#1}({#2})}
\newcommand{\lan}[2]{\mbox{\rm l.ann}_{#1}({#2})}

\begin{document}
\title{ \bf\normalsize ORE EXTENSIONS OVER DUO RINGS}
\author{ {\bf\normalsize Jerzy Matczuk\footnote{ The research was supported
by Polish KBN grant No. 1 P03A 032 27
}}\\
\normalsize Institute of Mathematics, Warsaw University,\\
 \normalsize Banacha 2, 02-097 Warsaw, Poland\\
 \normalsize E-mail: jmatczuk@mimuw.edu.pl}
\date{ }
\maketitle\markboth{ \bf  J.Matczuk}{ \bf ORE EXTENSIONS OVER DUO
RINGS}

\begin{abstract}
We show that there exist      noncommutative Ore extensions in
which every right  ideal is two-sided. This answers a problem
posed by Marks in \cite{Marks}. We  also  provide an easy
construction of one sided duo rings.

\end{abstract}

\section*{\bf \normalsize INTRODUCTION}
Hirano, Hong, Kim and Park proved in \cite{Hirano} that an
ordinary polynomial ring is one-sided duo only if it is
commutative.  Marks in \cite{Marks}, extended this result to Ore
extensions, by showing that if a noncommutative Ore extension
$\Tn$ which is  a duo ring on   one side exists,  then it has to
be right duo, $\si$  must be not injective and $\de\ne 0$  (see
Theorem \ref{Marks characterization of left duo}). He also
obtained a series of necessary conditions for the Ore extension to
be right duo (see Proposition \ref{necessary conditions}).

The aim of this paper is to show that noncommutative Ore
extensions which are right duo rings do exist and that the
necessary conditions obtained by Marks are not sufficient for the
Ore extension to be right duo.

In Section \ref{corner extensions which are duo} we investigate
corner extensions $R=A\oplus M$ of right duo rings $A$. In
particular, we show in Theorem \ref{construction of duo rings}
that $R$ is right duo provided the $(A-A)$-bimodule $M$ is simple
as a right $A$-module and faithful as a left $A$-module. This
theorem together with Lemma \ref{characterization of bimidules}
offer an easy way of constructing right duo rings which are not
left duo. Rings obtained in this way will serve us as coefficient
rings of Ore extensions which are right duo.

In Section \ref{skew derivations} we determine all
$\si$-derivations  for a suitably chosen corner extension
$R=A\oplus M$ and its endomorphism $\si$ (Theorem \ref{th
description of derivations}). This enable us to give, in
Proposition \ref{classification up to iso},  a classification of
Ore extensions $\Tn$ for $R$ and $\si$ as in Section \ref{corner
extensions which are duo}.

Finally, Section \ref{example} is devoted to description of Ore
extensions from the previous section which are right duo rings.

The construction from the paper gives not only Ore extensions
which are right duo. It provides also a ring R with an
endomorphism $\si$ such that, for any $n\in\mathbb{N}$, there
exists a $\si$-derivation $\de_n$ of $R$ such that   every right
ideal generated by a polynomial  of degree smaller than $n$ is a
two-sided ideal of $R[x;\si,\de_n]$ but there exists a polynomial
$f\in R[x;\si,\de_n]$ of degree $n$ for which $fR[x;\si,\de_n]$ is
not a left ideal. All of those examples
 satisfy necessary conditions, obtained by Marks, for an Ore
 extension to be right duo.

 \section{\bf \normalsize PRELIMINARIES}
 All rings considered in this paper are associative  with
 identity. Recall that a ring $R$ is called right (left) duo if every right
 (left) ideal of $R$ is a two-sided ideal.

For a subset $S$ of an $(R-R)$-bimodule $M$, $\lan{R}{S}$ will
stand for the left annihilator of $S$ in $R$, i.e.
$\lan{R}{S}=\{r\in R\mid rS=0 \}$. The right annihilator
$\an{R}{S}$ is defined similarly.

 An Ore extension of a ring $R$ is denoted by $\Tn$,
 where $\si$ is an endomorphism of $R$ and $\de$ is a $\si$-derivation,
 i.e. $\de \colon R\rightarrow R$ is an additive map such that
$\de(ab)=\si(a)\de(b)+\de(a)b$, for all $a,b\in R$.
 Recall that
elements of $\Tn$ are polynomials in $x$ with coefficients written
on the left. Multiplication in $\Tn$ is given by the
multiplication in $R$ and the condition $xa=\si(a)x+\de(a)$, for
all  $a\in R$.

We  say that a subset $S$ of  $R$ is $(\si-\de)$-stable, if
$\si(S)\subseteq S$ and $\de(S)\subseteq S$.

 For $a\in R$, the map
$d_a\colon R\rightarrow R$ defined by $d_a(r)=ar-\si(r)a$ is a
$\si$-derivation.  This $\si$-derivations is called  the inner
$\si$-derivation determined by the element $a$. A $\si$-derivation
$\de$ is called outer if it is not inner.

The following  fact is well-known (see for example Lemma
II.5.5\cite{Goodearl}):
\begin{lemma}\label{isomorphic extensions}
 Suppose that $\de_1$ and $\de_2$ are $\si$-derivations of a ring $R$.
 If $\de_1-\de_2$ is an inner $\si$-derivation, then the Ore extensions
 $R[x;\si,\de_1]$ and $R[x;\si,\de_2]$ are $R$-isomorphic.
\end{lemma}
The next two results come  from the paper \cite{Marks} of Marks.
\begin{theorem} \label{Marks characterization of left duo}
Suppose that one of the following conditions holds.
\begin{enumerate}
  \item The  Ore extension $R[x;\si,\de]$ is  a left duo ring;
 \item The  Ore extension $R[x;\si,\de]$ is  a right  duo ring
  and either $\si$ is injective  or $\de=0$.
\end{enumerate}
Then $R$ is commutative, $\si=\mbox{\rm id}_R$ and $\de=0$, i.e.
$R[x;\si,\de]=R[x]$ is a commutative polynomial ring.
\end{theorem}
\begin{prop}\label{necessary conditions}
 Suppose that $R[x;\si,\de]$ is a right duo ring which is noncommutative. Let
 $N=\bigcup_{i=1}^{\infty}\ker \si^i$.  Then:
\begin{enumerate}
\item $\de$ is an outer $\si$-derivation.
\item  $R$ is a right duo ring.
  \item  Every ideal of $R$ is $(\si-\de)$-stable.
  \item For any $r\in R$, the sequence $\{\si^n(r)\}_{n\in
  \mathbb{N}}$ is eventually constant.

  \item For any $r\in R$, the sequence $\{\si^n(\de(r))\}_{n\in
  \mathbb{N}}$ is eventually zero.
  \item  $0\ne N\subseteq J(R)$,
   where $J(R)$ denotes the Jacobson
  radical of $R$.
  \item The factor ring $R[x;\si,\de]/NR[x;\si,\de]$ is isomorphic
  to the commutative polynomial ring $(R/N)[x]$.

\end{enumerate}
\end{prop}
All statements but (1) from the above proposition come  from Lemma
7 and Theorem 11  of \cite{Marks}. The statement (1) is a direct
consequence of Lemma \ref{isomorphic extensions} and Theorem
\ref{Marks characterization of left duo}.

In \cite{Marks}, Marks also presented an  example of a ring  $R$
with an endomorphism $\si$ and $\si$-derivation $\de$ which
fulfils conditions $(3)\div (7)$. In this example  $\de$ is an
inner $\si$-derivation.

Let us observe that:
\begin{prop}\label{quasi-duo}
 Suppose that $R$, $\si$ and $\de$ possess  properties $(2)$, $(3)$
 and $(7)$ from Proposition \ref{necessary conditions}.
 If $N=\bigcup_{i=1}^{\infty}\ker \si^i$ is a nil ideal of $R$, then
 every maximal one-sided ideal of $\Tn$ is two-sided, i.e. $\Tn$
 is a quasi-duo ring.
\end{prop}
\begin{proof}
 Let $I$ be a nilpotent two-sided ideal  of $R$. By assumption,  $I$ is
 $(\si-\de)$-stable, so $I\Tn$ is also a nilpotent ideal of $\Tn$. In
 particular, $I\Tn$ is contained in the Jacobson
 radical $J$ of $\Tn$.

 Let $a\in N$. Since $R$ is a right duo
 ring and $a$ is a nilpotent element, $aR$ is a nilpotent two-sided ideal of
 $R$. Hence, by the above $N\Tn\subseteq J$ follows. This implies
 that $N\Tn$ is contained in any maximal one-sided ideal of $\Tn$.
Now, the thesis is an easy consequence of the fact that
$\Tn/(N\Tn)\simeq (R/N)[x]$ is a  commutative ring.
\end{proof}

Let $R^{\si}=\{r\in R\mid \si(r)=r\}$. Remark that statements
$(4)$ and $(5)$ of Proposition \ref{necessary conditions} say that
$R,\;\si$ and $\de$ are of very special form. Namely, in the
terminology  of Lam (Cf. Definition 2.15\cite{Lam2}), $R^{\si}$ is
a unital split corner ring of $R$, i.e. $R^{\si}$ is a unital
subring of $R$, $R=R^{\si}\oplus N$ as abelian groups and $N$ is
an ideal of $R$. The maps $\si$ and $\de$ satisfy: for any $r\in
N$, there exists $n\in\mathbb{N}$ such that $\si^n(r)=0$ and
$\de(R)\subseteq N$.

When seeking an  example of a noncommutative Ore extension $\Tn$
which is a right duo ring and  the coefficient ring $R$ is
one-sided noetherian,  one can restrict his attention to  more
specific rings:
\begin{prop}\label{nice corner}
 Let $R$ be either left or right noetherian ring.
 Suppose that $\Tn$ is a right duo ring which is noncommutative.
 Then there exists a  noncommutative Ore extension
 $R'[x;\si',\de']$, which is a right duo ring,
 such that:
\begin{enumerate}
  \item $R'=A\oplus M$ where $A$ is a unital split corner subring  of $R'$ with
 $M^2=0$ and $M\ne 0$.
  \item $\si'\colon R'\rightarrow R'$ is defined by $\si'(a+l)=a$,
  for any $a\in A$ and $l\in M$, and $\de'(R')\subseteq M$.
\end{enumerate}
\end{prop}
\begin{proof} We provide the proof in case $R$ is left noetherian.
The case $R$ right noetherian can be done using the same
arguments.

 We know, by Proposition \ref{necessary conditions},
  that $R^{\si}$ is a split corner subring of $R=R^{\si}\oplus
 N$, where $0\ne N=\bigcup_{i=1}^{\infty}\ker \si^i$ and $N\subseteq
 J(R)$.

Since $R$ is left noetherian, $N=\ker \si^m$ for some
$m\in\mathbb{N}$, say $m$ is the smallest such number.  If $m>1$,
then $\ker\si^{m-1}$ is a two-sided ideal of $R$ properly
contained in $N$. Proposition \ref{necessary conditions}(2)
implies that $\ker\si^{m-1}$ is a $(\si-\de)$-stable ideal.
Therefore $\si$
 and $\de$ induce an endomorphism $\si''$ and a $\si''$-derivation
 $\de''$ of $R/\ker\si^{m-1}$. The kernel of $\si''$ is equal to $N/\ker\si^{m-1}=N''\ne 0$.
  Notice that $\Tn/(\ker\si^{m-1}\Tn)\simeq (
 R/\ker\si^{m-1})[x;\si'',\de'']$ is a right duo
 ring as a homomorphic image of a right duo ring and it is not
 commutative, because $\ker \si''\ne 0$. Therefore, eventually
 replacing $R,\;\si$ and $\de$ by $R'',\;\si''$  and $\de''$,
 respectively, we may assume that $m=1$, i.e. $0\ne N=\ker \si\subseteq J(R)$.

Let us observe that $N^2\ne N$. Indeed, otherwise we would have
$N=N^2\subseteq J(R)N\subseteq N$, i.e. $J(R)N=N$. Notice that $N$
is finitely generated left $R$-module, as $R$ is left noetherian.
Thus, Nakayama's Lemma would imply $N=0$.

Now, one can make a similar reduction, as in the first part of the
proof,  using $N^2\ne N$ instead of $\ker\si^{m-1}$. This will
result in a noncommutative Ore extension $R'[x;\si',\de']$   such
that it is a right duo ring, $R'=A\oplus M$, where
$A=R^{\si'}=R^{\si}$, $M=\ker\si'$ and $M^2=0$. This proves (1)
and the first property of (2). The fact that $\de'(R')\subseteq M$
is a direct consequence of Proposition \ref{necessary
conditions}(5).
\end{proof}

In Section \ref{example} we will show that Ore extensions from the
above proposition do exist.

 Let us notice that it is easy to
construct rings described in the statement (2) from Proposition
\ref{nice corner}. In fact, if $A$ is a unital split corner
subring of $R=A\oplus M$ then $M$ is an $(A-A)$-bimodule.
Conversely, when we have a ring $A$ and an $(A-A)$-bimodule $M$,
then $A\oplus M$ is a ring with multiplication determined by
$M^2=0$, multiplication in $A$ and conditions $am=a\cdot m$,
$ma=m\cdot a$, for every $a\in A$ and $m\in M$, where $\cdot$
denotes the bimodule action of $A$ on $M$. This ring  can be
viewed also as a subring of a triangular ring $\begin{pmatrix}
  A & M \\
  0 & A
\end{pmatrix}$ consisting of all matrices of the form $\begin{pmatrix}
  a & m \\
  0 &a
\end{pmatrix}$, where $a\in A$ and $m\in M$.

Henceforth, while writing $R=A\oplus M$ we will always mean a ring
constructed as above, i.e. $A$ is a unital split corner subring of
$R=A\oplus M$ with  $M^2=0$. We will always assume that the
extension is nontrivial, i.e. $M\ne 0$. $\si$ will denote the
endomorphism of $R=A\oplus M$ defined by $\si|_A=\mbox{id}_A$ and
$\ker\si=M$. In the next two sections we will investigate when
 ring extensions of this kind are right duo  and we will describe all
 $\si$-derivations of such rings.

\section{ \bf \normalsize CORNER EXTENSIONS WHICH ARE RIGHT DUO RINGS}
\label{corner extensions which are duo}
 Throughout this section
$A$ is a unital split corner subring of $R=A\oplus M$, with
$M^2=0$.

We begin this section with the following  easy observation:
\begin{lemma}\label{reduction to left faithful}
 Suppose $R=A\oplus M$ is a right duo ring with $M\ne 0$. Then
 $I=\lan{A}{M}$ is a two-sided ideal of $R$ which is
 contained in $\an{A}{M}$.
 In particular, $M$ has a natural structure
 of $(A/I-A/I)$-bimodule and $R/I\simeq A/I\oplus M$.
\end{lemma}
\begin{proof} Let $I=\lan{A}{M}$.
Notice that $IR=I(A\oplus M)=IA\subseteq I$. Thus $I$ is a right
ideal of $R$ and, as $R$ is right duo, $I$ is a two-sided ideal.
Then $MI\subseteq M\cap I=0$, i.e. $I\subseteq \an{A}{M}$. Now it
is standard to complete the proof of the lemma.
\end{proof}
\begin{prop}\label{main prop}
Let $A$ be a right duo ring and $R=A\oplus M$, where  $M$ is
 an $(A-A)$-bimodule. Then:
\begin{enumerate}
  \item Suppose that $a\in A $ is such that $aM=M$. Then, for any
  $m\in M$, $(a+m)R=aA+M$
  is a two-sided ideal of $R$.
  \item Suppose that any right $A$-submodule of $M$ is a $(A-A)$-subbimodule of $M$
  and $aM=M$, for any $0\ne a\in A$. Then $R=A\oplus M$ is a right duo
  ring.
 \end{enumerate}
\end{prop}
\begin{proof}
 (1). Let  $a\in A$ be such that $aM=M$. Then, for any
  $m\in M$, $(a+m)M=aM=M$.
 Therefore $M\subseteq (a+m)R$. This yields that  also $aA\subseteq
 (a+m)R$ and we have $aA+M\subseteq (a+m)R\subseteq aR+mR\subseteq aR+M=aA+M$.
 This shows that $(a+m)R=aA+M$. Since $A$ is a right duo ring,
 $aA$ is a two-sided ideal of $A$. Now it is clear that $(a+m)R$ is
 a two-sided ideal of $R$.

 (2). The assumptions imposed on $M$ imply that, for any $m\in M$,
 $mR$ is a two-sided ideal of $R$ while, by the statement (1),
 $(a+m)R$ is a two-sided ideal of $R$ for any $0\ne a\in A$. This
 gives (2).
 \end{proof}
 \begin{remark} Let $R=A\oplus M$ be a ring satisfying the
 assumptions of the statement $(2)$ from the above proposition.
  Then,  using the  proposition, one can easily give a
  description of the lattice of two-sided ideals of the righ duo ring $R$
  in terms of the lattice of two-sided ideals of $A$ and the
  lattice of $(A-A)$-subbimodules of $M$.
 \end{remark}
 As an  application of  Proposition \ref{main prop} we obtain the following:
\begin{theorem}\label{construction of duo rings}
 Let $A$ be a right duo ring and $R=A\oplus M$, where  $M$ is
 an $(A-A)$-bimodule such that
  $M$ is faithful as a
  left $A$-module and simple as a right $A$-module. Then:
\begin{enumerate}
  \item[1.]    $R$  is a   right  duo ring.
  \end{enumerate}
  Moreover:
\begin{enumerate}
    \item[2.] $R$ is left duo iff  $M$ is  faithful as a
  right $A$-module (i.e. $A$ is  a division ring) and simple as a
  left $A$-module.
  \end{enumerate}
  \end{theorem}
\begin{proof}(1).
 Let $0\ne  a\in A$.  Since the module
 $_AM$ is faithful, $aM$ is a nonzero submodule of the simple
 right $A$-module $M_A$, so $aM=M$ follows. Now, it is easy to see that
  the thesis is a consequence of Proposition \ref{main prop}(2).

(2). Suppose that $R$ is left duo. Then, by the left hand version
of Lemma \ref{reduction to left faithful}, $I=\an{A}{M}\subseteq
\lan{A}{M}=0$.  Hence $I=0$ and $M$ is faithful as a right
$A$-module. Since $M_A$ is also simple, $A$ is a right
 primitive, right duo ring. Thus $A$ is a division ring.

 If $0\ne N\subseteq {_A}M$ is an
 $A$-submodule of $_AM$ then $N$ is a two-sided ideal of $R$ since $R$
 is left duo. Thus, in particular,  $N$ is also a submodule of the
 simple $A$-module $M_A$ and $N=M$ follows, showing that the left  $A$-module $_AM$ is
   simple.

 On the other hand, if $M$ is simple and faithful both as a left and
 right $A$-module, then $A$ is a division ring, as $A$ is right duo, and $M$ is the only
 proper one-sided ideal  of $R$. In particular, $R$ is left duo.
\end{proof}
The  examples of right duo rings  given by  the above theorem are,
in some sense, minimal. Lemma \ref{reduction to left faithful}
shows that a right duo ring $R=A\oplus M$ has a quotient  of the
form $A'\oplus M'$, where $_{A'}M$ is a faithful $A'$-module. The
following proposition offers another reduction.
\begin{prop}\label{reduction}
Suppose that $R=A\oplus M$ is a right duo ring and the
(A-A)-bimodule $M$
 has a maximal $(A-A)$-subbimodule.  Then there exists an
ideal $J$ of $R$ such that $R/J\simeq A'\oplus M'$, for some right
duo ring $A'$, where $M'$ is an $(A'-A')$-bimodule which is simple
as a right $A'$-module and faithful as a left $A'$-module.
\end{prop}
\begin{proof}
 Let $N$ be a maximal (A-A)-subbimodule of $M$. Then $N$ is a
 two-sided ideal of $R$ and $R/N\simeq A\oplus M'$, where $M'$ denotes the
 quotient $(A-A)$-bimodule $M/N$.

 Notice that if $0\ne W_A$ is a submodule of $M'_A$, then $W_A$ is a
 two-sided ideal of $R/N$ as it is a right ideal of a right duo ring $R/N$.
  In particular, $AW_A\subseteq W_A$, i.e. $0\ne W_A$ is
 subbimodule of a simple bimodule $M'$. This means that $W=M'$ and shows that $M'$ is
 simple as right $A$-module. Now the proposition is a direct consequence
 of Lemma \ref{reduction to left faithful}. Remark that the ideal
 $J$ from the proposition is equal to $\lan{A}{_A(M/N)}+N$.
\end{proof}
In order to be able to make use of Theorem \ref{construction of
duo rings}, we need the following:
\begin{lemma}\label{characterization of bimidules}
 Let $A$ be a right duo ring. The following conditions are
 equivalent:
\begin{enumerate}
  \item  There exists an $(A,A)$-bimodule $M$ such that $M$
  is faithful as left $A$-module and simple as right $A$-module.
  \item  There exist a right primitive ideal $P$ of $A$ and an
  injective homomorphism $\phi\colon A\rightarrow A/P$.
\end{enumerate}
\end{lemma}
\begin{proof}
 $(1)\Rightarrow (2)$. Let $P$ denote the annihilator of $M_A$.
 Then $P$ is right primitive ideal of $A$ and $A/P$ is a division ring as $A$
 is a right duo ring. This means that, for any $0\ne m\in M$,
 $\an{A}{m}=P$.

 Let us fix $0\ne m\in M$ and consider $M$ as
 $(A-A/P)$-bimodule. Then $M=m(A/P)$ and for any $a\in A$,
 $am=m\phi_m(a)$ for a suitable element $\phi_m(a)\in A/P$.
  Notice that, because $\an{A/P}{m}=0$,
  the element $\phi_m(a)$ is uniquely determined by $a$. Thus we
  have a well-defined  map $\phi=\phi_m\colon A\rightarrow A/P$. It is standard to
  check that $\phi$ is a ring homomorphism. If $\phi(a)=0$, then
  $0=m\phi(a)(A/P)=am(A/P)=aM$. Hence $a=0$ follows, as the left $A$-module
  $_AM$ is faithful. This shows that $\phi$ is injective.

$(2)\Rightarrow (1)$. Let $P$ be a right primitive ideal of $R$
and $\phi\colon A\rightarrow A/P$ an injective homomorphism. Then,
as $A$ is right duo, $A/P$ is a division ring. Let $M$ be the one
dimensional right vector space over $A/P$. Let us fix $0\ne m\in
M$ and define left $A$ module structure on $M$ by setting
$a\cdot(mr)=m\phi(a)r$, for any $a\in A$ and $r\in A/P$. It is
standard to check that this determines an $(A-A/P)$-bimodule
structure on $M$.  Notice that if $aM=0$, then $m\phi(a)=0$ and
$a=0$ follows, as $\phi$ is injective and $\an{A/P}{m}=0$. This
induces the desired $(A-A)$-bimodule structure on $M$.
\end{proof}
\begin{remark}\begin{enumerate}
 \item In the proof of the implication $(1)\Rightarrow (2)$ from the
above lemma, different choices of the   element $0\ne m\in M$ give
different homomorphisms $\phi_m\colon A\rightarrow A/P$. In fact,
one can check that $\phi_{mr}=r^{-1}\phi_mr$, for any $0\ne r\in
A/P$.
\item The equivalence from the above lemma holds under a slightly
weaker assumption that every maximal right ideal of $A$ is
two-sided, i.e. $A$ is a right quasi-duo ring.
\end{enumerate}
\end{remark}
\begin{corollary}
 Suppose that a right duo ring $A$ satisfies one of the equivalent
 conditions of the above lemma. Then:
\begin{enumerate}
  \item $A$ is a right Ore domain.
  \item If $A$ is an algebra over a field $K$ then  $A$ is a division
  $K$-algebra, provided $A$ satisfies
  Nullstellensatz, i.e.  for any simple $A$-module $N_A$, the division $K$-algebra
End$_A(N_A)$ is algebraic over $K$ (for example when $A$ is a
  finitely generated commutative $K$-algebra or $\dim_KA< \# (K)-1$,
  as cardinal numbers).

\end{enumerate}
\end{corollary}
\begin{proof}
 The statement (1) is clear as $A\simeq\phi(A)$ is a subring of
 the division ring $A/P$ and a domain which is a right duo ring
  is always a right Ore domain.

(2) Suppose the $K$-algebra $A$ satisfies Nullstellensatz and
 let $M$ denote the
$(A-A)$-bimodule from Lemma \ref{characterization of bimidules}.
Then, since $M_A$ is simple, End$_A(M_A)$ is an algebraic division
algebra over $K$.

 Notice that $\psi\colon A\rightarrow \mbox{\rm
End}_A(M_A)$ defined by $\psi(a)(m)=am$, for $a\in A$ and $m\in
M$, is a $K$-algebra homomorphism. Moreover $\psi$ is injective,
since $_AM$ is a faithful as a left $A$-module. Therefore $A\simeq
\psi(A)$ is a domain which is algebraic over $K$. This means that
$A$ is a division algebra.
\end{proof}
Of course, in general $A$ does not have to be a division ring if A
possesses  an $(A-A)$-bimodule which is simple as a right
$A$-module and faithful as left $A$-module.

\begin{example}\label{example 1}
 Let $A=K[x]$, where   $K=F(X)$ denotes the field
 of rational functions over a field $F$ in the set $X=\{x_i\mid i=0,1,\ldots\}$ of
 indeterminates.
  Then the $F$-linear homomorphism  $\phi\colon A\rightarrow K$
 defined by $\phi(x)=x_1$ and $\phi(x_i)=x_{i+2}$, for $i=0,1,\ldots$, is injective.
 Thus, by Lemma \ref{characterization of bimidules}
  and Theorem \ref{construction of duo rings}, $R=A\oplus M$ is a right duo ring which
  is not left duo, where $M=K$ has the
  $(A-A)$-bimodule structure given by $a\cdot
  m=\phi(a)m$ and $m\cdot a=m\overline{a}$,
  for $m\in M$ and $a\in A$, where $\ov{a}$ denotes the
  canonical image of $a$ in $A/(x)=K$.
\end{example}

Theorem \ref{construction of duo rings} together with Lemma
\ref{characterization of bimidules} offers an easy way of
constructing right duo rings which are not left duo.

Let us notice,  that when $K=F(x)$ is the field of rational
functions in indeterminate $x$ and $\phi\colon K\rightarrow K$ is
an $F$-homomorphism given by $\phi(x)=x^2$, then the resulting
right duo ring is exactly an old  example of Asano (see Exercise
22.4A\cite{Lam1} and comments hereafter).

\section{\bf \normalsize SKEW DERIVATIONS OF $A\oplus
M$}\label{skew derivations}
  Henceforth $A$ will stand for a commutative domain, $P$
for a maximal ideal of $A$, $K$ will denote the field $A/P$ and
$\phi\colon A\rightarrow K$ a fixed injective homomorphism of
rings. For any element $a\in A$, $\ov{a}$ will denote   the
canonical image of $a$ in $K=A/P$.

By Lemma \ref{characterization of bimidules}, the right $K$ vector
space $vK$ with the basis $\{v\}$ has a structure of
$(A-A)$-bimodule given by $a\cdot vk=v\phi(a)k$ and $vk\cdot
a=vk\ov{a}$, for any $a\in A$ and $k\in K$. $vK$ is faithful as a
left $A$-module and simple as a right $A$-module. Thus, by Theorem
\ref{construction of duo rings}, $R=A\oplus vK$ is a right duo
ring.

From now on,  $\si :R\rightarrow R$ stands for the  endomorphism
of $R$ given by $\si(a+v\ov{l})=a$, for any $a, l\in A$.
\begin{lemma}\label{description of inner derivations}
 Let $d_y$ denote the inner $\si$-derivation of $R$ determined by
 the element $y=c+v\ov{m}\in A\oplus vK=R$, where $c,m\in A$. Then:
\begin{enumerate}
  \item $d_y(a+v\ov{l})= v\phi(c)\ov{l}+ v\ov{m}(\ov{a}-\phi(a))\in
 vK$, for any $a+v\ov{l}\in R$. In particular, $d_y(v)=v\phi(c)$.
  \item If $d_y(A)=0$, then $d_y(a+v\ov{l})=v\phi(c)\ov{l}$, for any $a+v\ov{l}\in R$.

\end{enumerate}
\end{lemma}
\begin{proof}
(1) By definition,
\begin{equation*}
\begin{split}
  d_y(a+v\ov{l}) & =y(a+v\ov{l})-\si(a+v\ov{l})y=
 ca+v\ov{m}\ov{a}+v\phi(c)\ov{l}-ac-v\phi(a)\ov{m}\\
  & =v\phi(c)\ov{l}+ v\ov{m}(\ov{a}-\phi(a)).
\end{split}
\end{equation*}

Taking $a=0$ and
 $l=1$, we obtain $d_y(v)=v\phi(c)$.

 (2) If $d_y(A)=0$ then, by the statement (1),
 $v\ov{m}(\ov{a}-\phi(a))=0$, for all $a\in A$, i.e.
 $d_y(a+v\ov{l})=v\phi(c)\ov{l}$, for $a,l\in A$.
\end{proof}
For $\vm\in K$, define $\dw\colon R\rightarrow R$ by setting
 $\dw(a+v\ov{l})=v\vm\ov{l}$, for any $a,l\in A$. Keeping the
 above notion, we have:
\begin{lemma}\label{outer derivations}
For any $\vm\in K$, $\dw$ is a $\si$-derivation of $R=A\oplus vK$.
Moreover $\dw$ is an inner $\si$-derivation iff  $\vm\in \phi(A)$.
\end{lemma}
\begin{proof} Let $r_1=a+v\ov{l}, r_2=b+v\ov{m}\in R=A\oplus vK$, where $a,b,l,m\in
A$.
 By making direct computations we have:
$$\dw(r_1r_2)=\dw(ab+v(\phi(a)\ov{m}
+\ov{l}\ov{b}))=v\vm(\phi(a)\ov{m}+\ov{l}\ov{b})\;\;\;
\mbox{and}$$
 $$\si(r_1)\dw(r_2)+\dw(r_1)r_2=a\cdot
v\vm\ov{m}+v\vm\ov{l}(b+v\ov{m})=v\vm\phi(a)\ov{m}+v\vm\ov{l}\ov{b}.$$
This shows that $\dw$ is a $\si$-derivation.

Notice that $\dw(A)=0$ and, by Lemma \ref{description of inner
derivations}(2), every inner $\si$-derivation of $R$ such that
$\de(A)=0$  is  of the form $\de(a+v\ov{l})= v\phi(c)\ov{l}$ for
suitable $c\in A$. Hence $\dw$ is inner iff there is $c\in A$ such
that $\vm=\phi(c)$. This completes the proof of the lemma.
\end{proof}

Recall that Lemma \ref{description of inner derivations} describes
all inner $\si$-derivations of $R$. Therefore, the following
theorem gives a description of all $\si$-derivations of $R$.
\begin{theorem} \label{th description of derivations}
Let $\de$ be a nonzero $\si$-derivation of $R=A\oplus vK$. Then:
\begin{enumerate}
  \item There exists $\vm\in K$ such that $\de(v)=v\vm$.
  \item If $\de(vK)=0$, then  one of the following
  conditions holds:
  \begin{enumerate}
   \item $\de$ is an inner $\si$-derivation of $R$
\item $\phi=\mbox{\rm id}_K$, i.e. $R=K\oplus vK$ is a
  commutative ring, $\de$ is  an outer $\si$-derivation
   and there exists a derivation $d$ of the field $K$ such that
  $\de(a+vb)=vd(a)$, for any $a,b\in K$.
    \end{enumerate}

  \item  Let $\vm \in K$ be such that $\de(v)=v\vm$. Then
  $(\de-\dw)(vK)=0$, i.e. $\de-\dw$ is a $\si$-derivation
  satisfying the assumption of the statement (2).

 \end{enumerate}
\end{theorem}
\begin{proof}
(1). Let $\de(v)=a+v\vm$, where $a,\om\in A$. Since $\si(v)=0$, we
have: $0=\de(0)=\de(v^2)=\si(v)\de(v)+\de(v)v=\de(v)v=v\phi(a)$.
Hence $a=0$ follows, as $\phi$ is injective.  This gives (1).

(2). Let $\de$ be a  $\si$-derivation of $R=A\oplus vK$ such that
$\de(vK)=0$. First we claim that $\de(A)\subseteq vK$. To this
end, let $a\in A$ and $\de(a)=c+v\ov{s}\in R$, for some $c,s\in
A$. Then $0=\de(v\phi(a))=\de(a\cdot v)=
\si(a)\cdot\de(v)+\de(a)\cdot v= (c+v\ov{s})v= v\phi(c)$. Hence,
as $\phi$ is injective, $c=0$. This proves the claim.

Since  $\de\ne 0$ and $\de(vK)=0$, we may pick an element $a_0\in
A$ such that $ \de(a_0)\ne 0$. By the first part of the proof
$\de(a_0)=v\ov{s}$, for some $s\in A$.

Let $b\in A$. Then $\de(b)\in vK$, so $a_0\cdot
\de(b)=\de(b)\phi(a_0)$ and $\de(b)\cdot a_0=\de(b)\ov{a_0}$.
Computing $\de(a_0b)=\de(ba_0)$ we obtain
$v\ov{s}\ov{b}+a_0\cdot\de(b)=\de(b)\cdot a_0+v\phi(b)\ov{s}$.
Using this, one can see that
\begin{equation}\label{1}
v\ov{s}(\ov{b}-\phi(b))=\de(b)(\ov{a_0}-\phi(a_0))\;\; \mbox{\rm
for any } b\in A.
\end{equation}
Let $c\in A$ be such that $\ov{c}=\ov{a_0}-\phi(a_0)\in K$.

If $\ov{c}=0$  then, using the equation (\ref{1}) and the fact
that $\ov{s}\ne0$, we get $\ov{b}-\phi(b)=0$, for any $b\in A$.
Recall that  $\phi\colon A\rightarrow A/P=K$ is injective, and
$\ov{b}$ is a
 natural image of $b\in A$ in $K$. Therefore, $P$ has to be equal
 to $0$,  i.e. $A=K$ is a field  and $\phi=\mbox{id}_K$. Then
  $R=K\oplus vK$ is a commutative ring (isomorphic to
  $K[x]/(x^2)$). By the first part of the proof of (2), $\de(K)\subseteq
  vK$. Therefore, for any $a\in K$,  $\de(a)=vd(a)$ for a suitable
  element $d(a)\in K$. It is clear that the element $d(a)$ is
  uniquely determined by $a$. Moreover, as also $\de(vK)=0$,
  $\de(a+vb)=d(a)$, for any $a,b\in K$. Now, it is standard to
  check that the map $\de\colon R\rightarrow R$, defined by the above formula, is a
  $\si$-derivation of $R$ iff $d\colon K\rightarrow K$ is a
  derivation of the field $K$. Moreover, by Lemma
  \ref{description of inner derivations}(1), such nonzero
  $\si$-derivation is always outer,
   i.e. the statement (b) holds.

 Suppose $0\ne\ov{c}\in K$. Because  $K$ is a field, there exists $c'\in
 A$ such that $\ov{c}\ov{c'}=1$. Then, again making use of the
 equation (\ref{1}), we obtain
 $\de(b)=v\ov{s}\ov{c'}(\ov{b}-\phi(b))$ for any $b\in B$. Now,
 the fact that $\de(vK)=0$ and Lemma \ref{description of inner
 derivations}(1) yield that $\de=d_{y}$ is the inner $\si$-derivation
 determined by the element $y=v\ov{s}\ov{c'}\in R$, i.e. the statement (a) holds.
  This completes the proof of (2).

(3). By (1), there is $\vm\in K$ such that $\de(v)=v\vm$. Notice
that $(\de -\dw)(v)=0$ and $\si(v)=0$. Therefore $$(\de
-\dw)(v\ov{l})=(\de -\dw)(v\cdot l)=\si(v)\cdot (\de -\dw)(l)+
(\de -\dw)(v)\cdot l=0,$$ for any $l\in A$. This means that $(\de
-\dw)(vK)=(\de -\dw)(v\cdot A)=0$ and the statement (3) follows.
\end{proof}
The above theorem together with Lemma \ref{outer derivations} give
us immediately the following:
\begin{corollary}\label{description of outer der}
 Suppose that $R=A\oplus vK$ is noncommutative. For a $\si$-derivation
 $\de$ of $R$, the following conditions are equivalent:
\begin{enumerate}
  \item $\de$ is an outer $\si$-derivation of $R$.
  \item There exist $\vm\in K$ and $y\in R$, such that:\\
  (i) $\vm\not \in \phi(A)$;\\
  (ii) $\de=\dw+d_y$.
\end{enumerate}
\end{corollary}

If $R$ is commutative (i.e. $\phi=\mbox{id}_A$) then, By Lemma
\ref{outer derivations}, $\de_{\om}$ is an inner derivation of
$R$, for any $\vm\in K$. Thus, by the above theorem, we also get:
\begin{corollary}
 Suppose $R=A\oplus vK$ is commutative. Then $R=K\oplus vK$ and
 for a $\si$-derivation
 $\de$ of $R$, the following conditions are equivalent:
 \begin{enumerate}
  \item $\de$ is an outer $\si$-derivation of $R$.
  \item There exist a nonzero derivation $d$ of $K$, such that
  $\de(a+vb)=vd(a)$, for any $a,b\in K$.
\end{enumerate}
\end{corollary} As a direct application of the above Corollaries and Lemma
\ref{isomorphic extensions} we obtain   the following
classification of Ore extensions over our ring $R=A\oplus vK$:

\begin{prop}\label{classification up to iso} Let $\de$ be a
$\si$-derivation of $R=A\oplus vK$. Then:
\begin{enumerate}
  \item  Suppose that $R$ is
noncommutative.
 Then the Ore extension $\Tn$ is  $R$-isomorphic  either to $R[x;\si]$
 or to $R[x;\si,\dw]$, for some $\vm\in K\setminus \phi(A)$,
 where $\dw(a+v\ov{l})=v\vm\ov{l}$, for any $a,l\in A$.
  \item Suppose that $R$ is commutative. Then $R=K\oplus vK$ and
  the Ore extension
   $\Tn$ is  $R$-isomorphic  either to $R[x;\si]$ or to
   $R[x;\si,\hat{\de}]$, where $\hat{\de}(a+vb)=vd(a)$, for any $a,b\in
   K$, and $d$ denotes a nonzero derivation of $K$.
 \end{enumerate}

\end{prop}

\section{\bf \normalsize ORE EXTENSIONS WHICH ARE RIGHT
DUO RINGS}\label{example}  We will continue to use the notation
from the previous section.
 Recall that  $A$ is a commutative domain, $K=A/P$, where
 $P$ is a fixed maximal ideal of $A$. $R=A\oplus vK$ is a split
 corner
extension  of $A$  with $(vK)^2=0$. The left and right actions of
$A$ on the right $K$-linear vector space $vK$ are given by $a\cdot
v=v\phi(a)$ and $v\cdot a=v\ov{a}$, for $a\in A$.

Recall that $\si :R\rightarrow R$ denotes the endomorphism of  $R$
given by $\si(a+v\ov{l})=a$ for any $a,l\in A$. For a fixed
element $\vm\in K$, $\dw$ stands for the $\si$-derivation of $R$
defined in Lemma \ref{outer derivations}, i.e.
$\dw(a+v\ov{l})=v\vm\ov{l}$, for any $a,l\in K$.

\begin{lemma}\label{bbasic} For any $a\in A$ and $f\in \T$, we have:
\begin{enumerate}
\item $xa=ax$ and $xv=v\vm$;
\item $vfv=0$
\end{enumerate}
In particular, $v\T$ is a two-sided ideal of $\T$ with
$(v\T)^2=0$.
\end{lemma}
\begin{proof}
 The easy proof is left to the reader.
\end{proof}
For $\vm\in K$, let $\p$ denote the extension of $\phi$ to a
homomorphism $\p\colon A[x]\rightarrow K\subseteq K[x]$ given by
$\p(x)=\vm$.

We also extend $\bar{} \colon A\rightarrow K$ to a homo\-morphism
$\bar{}\colon A[x]\rightarrow K[x]$, by setting $\ov{x}=x$.

Let $vK[x]$ denote the free right $K[x]$-module generated by the
element $v$. Then $vK[x]$ has a  structure of an
$(A[x]-A[x])$-bimodule given by $f\cdot v=v\p(f)$ and $v\cdot
f=v\ov{f}$, for  $f\in A[x]$. Thus we can consider the ring
$T_{\om}=A[x]\oplus vK[x]$. With the help of   Lemma \ref{bbasic},
one can easily check that:
\begin{lemma}\label{isomorphism}
\begin{enumerate}
  \item  If $f\in A[x]\subseteq \T$, then $fv=v\phi_{\om}(f)$
  \item The map $\Phi\colon \T\rightarrow T_{\om}$ defined by
$\Phi((a+v\ov{l})x^k)=ax^k\oplus v\ov{l}x^k$, for any $a,l\in A$,
is an isomorphism of rings.
\item If $N$ is an $A[x]$-submodule of the right $A[x]$-module
$vK[x]\subseteq T_{\om}$, then $N$ is a $(A[x]-A[x])$-subbimodule
of $vK[x]$.
  \end{enumerate}
\end{lemma}

\begin{definition}
 For any polynomial
$f=\sum_{k=0}^{n}(a_k+v\ov{l_k})x^k\in\T$ we set:
\begin{enumerate}
  \item  $f_A=\sum_{k=0}^{n}a_kx^k\in A[x]\subseteq \T$ and $f_v=f-f_A$.
  \item $D_f=\sum_{k=0}^{n}\phi(a_k)\vm^k\in K$, that is
  $D_f=\phi_{\om}(f_A)$.
 \end{enumerate}
\end{definition}

Notice that if the element $\vm\in K$ is transcendental over the
subfield generated by $\phi(A)\subseteq K$, then $D_f=0$ iff
$a_k=0$ for all $0\leq k\leq n$, i.e. $f=f_v\in v\T$.

Combining Lemmas \ref{bbasic} and \ref{isomorphism}(1) we get the
following:
\begin{remark} \label{formula} Let
 $f\in \T$. Then:
 $fv=f_Av=vD_f$. In particular, $fv=0 $ iff $ D_f=0$.
 \end{remark}

Now we are in position  to prove the following:
\begin{prop}\label{description of ideals}
 For a polynomial $f\in \T$, the following conditions are
 equivalent:
\begin{enumerate}
  \item  $f\T$ is a two-sided ideal of $\T$;
  \item One of the following conditions holds:
  \begin{enumerate}
  \item $D_f\ne 0$;
  \item $f_A=0$, i.e. $f\in v\T$;
  \item $vf=0$ and $f=f_A$.
  \end{enumerate}
\end{enumerate}
\end{prop}
\begin{proof}Let $f=\sum_{k=0}^{n}(a_k+v\ov{l_k})x^k\in\T$,
 with $a_n+v\ov{l_n}\ne 0$.

$(2)\Rightarrow (1)$. By Lemma \ref{isomorphism}, the ring $\T$ is
isomorphic to $T_{\om}=A[x]\oplus vK[x]$.  Let
$f_A=\sum_{k=0}^{n}a_kx^k\in A[x]$. Notice that
$D_f=D_{f_A}=\phi_{\om}(f_A)$.

Suppose that $D_f\ne 0$. Then, by Lemma \ref{isomorphism},
$f_AvK[x]=v\phi_{\om}(f_A)K[x]=vK[x]$.  Therefore, by Proposition
\ref{main prop}(1), $f_AT_{\om}=fT_{\om}$ is a two-sided ideal of
$T_{\om}$. This yields that $f\T$ is a two-sided ideal of $\T$,
provided $D_f\ne 0$.

Suppose $f_A=0$, i.e. all coefficients of $f$ are in $vK$. In this
case,  Lemma \ref{isomorphism}(3) implies that $f\T$ is a
two-sided ideal of $\T$.

Finally, suppose that the statement (c) of (2) holds. Let
$g\in\T$. Since $vf=0$ and $f=f_A$, we obtain
$gf=g_Af=g_Af_A=f_Ag_A=fg_A$. This shows that $f\T$ is a two-sided
ideal of $\T$.

$(1) \Rightarrow (2).$ Let $f R[x; \sigma, \delta_{\omega}]$ be a
two-sided ideal of $R[x; \sigma, \delta_{\omega}],$ where $f =
\sum_{k=0}^{n} (a_k + v \overline{l_k}) x^k \in R[x; \sigma,
\delta_{\omega}]$ with $a_n + v \overline{l_n} \neq 0.$ Suppose
that $D_f = 0$.  We shall prove that either (b) or (c) of the
statement (2) holds.

By Remark \ref{formula}, for any $g \in R[x; \sigma,
\delta_{\omega}]$ we have $fg = f g_A.$ We claim that $vf = 0.$
Since $f R[x; \sigma, \delta_{\omega}]$ is a two-sided ideal, $vf
\in f R[x; \sigma, \delta_{\omega}].$ Thus, for some $g \in R[x;
\sigma, \delta_{\omega}]$ we have:
\begin{equation}
vf = v f_A = fg = f g_A. \label{ClaimEquation1}
\end{equation}
Then $vf = (f_A + f_v)g_A = f_A g_A + f_v g_A$ implies $0 = f_A
g_A$ in the domain $A[x],$ whence $f_A = 0$ or $g_A = 0.$ In
either case, Equation (\ref{ClaimEquation1}) implies $vf = 0.$

It now suffices to show that either $f_v = 0$ or $f_A = 0.$ To
this end, assume that $f_A \ne 0.$ Since $vf = 0,$ we have $f_A
\in P[x] \subseteq A[x]$ (where $A/P=K$). Thus, in particular, $P
\ne 0$. Choose a nonzero element $p \in P.$ Since $pf \in f R[x;
\sigma, \delta_{\omega}],$ for some $h \in R[x; \sigma,
\delta_{\omega}]$ we have:
\begin{equation}
pf = f_A p + f_v \phi(p) = fh = f h_A = f_A h_A + f_v h_A.
\label{ClaimEquation2}
\end{equation}
Hence $f_A p = f_A h_A$ in the domain $A[x],$ so $h_A = p\in P$.
Therefore $f_vh_A=f_vp=f_v\ov{p}=0$ and Equation
(\ref{ClaimEquation2}) implies $f_v \phi(p) = 0$. Since $p \ne 0$
and $\phi$ is injective, $f_v = 0$ follows. This completes the
proof of the proposition.
\end{proof}

Recall that for $f=\sum_{k=0}^{n}(a_k+v\ov{l_k})x^k\in\T$,
$D_f=\sum_{k=0}^{n}\phi(a_k)\vm^k$. Notice that if either
$\ov{\om}$ is transcendental  over the subfield
$\widehat{\phi(A)}$ of $K$ generated by $\phi(A)$ or $\vm$ is
algebraic over $\widehat{\phi(A)}$ of degree greater than $\deg
f$, then  $D_f= 0$ iff $f\in v\T$. Hence, by the above
proposition, $f\T$ is a two-sided ideal of $\T$.

 On the other
hand, if $\ov{\om}$ is algebraic, say of degree $n$, then there
exists a polynomial $g=\sum_{k=0}^na_kx^k\in A[x]\subseteq \T$ of
degree $n$, such that $D_g=0$. Then, by the above proposition, the
right ideal  $(g+v)\T$ is not two-sided. Therefore we obtain:
\begin{corollary}\label{when ore ext is duo}Let $\vm \in K$ and $\widehat{\phi(A)}$
denote the subfield of $K$ generated by $\phi(A)$. Then:
\begin{enumerate}
  \item  If $\vm$ is transcendental over $\widehat{\phi(A)}$
   then $f\T$ is a two-sided ideal of $\T$, for any $f\in\T$, i.e.
  $\T$ is a right duo ring.
  \item If $\vm$ is algebraic of degree $n+1$ over $\widehat{\phi(A)}$, for some $n\geq
  0$, then:
\begin{enumerate}
  \item for every polynomial $f\in \T$ of degree $\deg(f)\leq n$, $f\T$ is a two-sided
  ideal of $\T$;
  \item  there exists a polynomial $f\in \T$ of degree $n+1$ such
  that $f\T$ is not a two-sided ideal of $\T$.
\end{enumerate}
\end{enumerate}
\end{corollary}

When $P=0$, i.e. $R=K\oplus vK$, then $vf\ne 0$, for any
polynomial $f\in\T$ with $f_A\ne 0$. Thus, in this case,
Proposition \ref{description of ideals} boils down to:
\begin{corollary}\label{case A is a field}
Suppose $R=K\oplus vK$ and $f\in \T$. Then $f\T$ is a two-sided
ideal of $\T$ iff either $D_f\ne 0$ or $f\in v\T$.
\end{corollary}
 Now we are in
position to prove the following:
\begin{theorem}\label{theorem duo rings}
 Let $A$ be a commutative domain with a maximal ideal $P$, $\phi\colon A\rightarrow
 A/P=K$ an injective homomorphism and $R=A\oplus vK$ the
 associated unital split corner extension of $A$. Then:
 \begin{enumerate}
 \item  $\Tn$ is a quasi-duo ring.
 \item The following statements are equivalent:

 \begin{enumerate}
  \item $\Tn$ is a right duo ring;
  \item There exists  $\vm \in K$, such that $\vm$ is transcendental over the subfield of
 $K$ generated by $\phi(A)$ and $\Tn$ is $R$-isomorphic to $\T$.
 \end{enumerate}

 \end{enumerate}
\end{theorem}
\begin{proof}
(1). By virtue of Proposition \ref{classification up to iso}, it
is enough to prove the statement in   case the Ore extension $\Tn$
is as described in the proposition. Then, in any case,  $I=vK$ is
a nilpotent $(\si-\de)$-stable ideal of $R$ and  $\Tn/I\Tn\simeq
A[x]$ is commutative. Now, one can complete the proof of (1) using
similar arguments as in the proof of Proposition \ref{quasi-duo}.

(2). The implication $(b)\Rightarrow (a)$ is given by Corollary
\ref{when ore ext is duo}(1).

 $(a)\Rightarrow (b)$. By Theorem \ref{Marks characterization of left
duo}(2), the Ore extension $R[x;\si]$ is never right duo. Thus, in
view of Proposition \ref{classification up to iso} and Corollary
\ref{when ore ext is duo}(2), it is enough to show that if $R$ is
commutative and $\de$ is of the form described in Proposition
\ref{classification up to iso}(2), then $\Tn$ is not right duo. To
this end, suppose that $R=K\oplus vK$ is commutative, $d$ is a
derivation of $K$ and $\de(a+vb)=vd(a)$, for any $a,b\in K$. We
claim that $vx\not \in x\Tn$, i.e. $x\Tn$ is not a left ideal.
Suppose that $vx=xg$, for some polynomial $g\in\Tn$. Since
$\si(v)=\de(v)=0$, $xv=0$ and we may assume that
$g=\sum_{i=0}^na_ix^i$, where $0\ne a_n,a_{n-1},\ldots , a_0\in
K$. Then $\deg (vx)=\deg (xg)=n+1$ and $g=a\in K$ follows. Hence
$vx=ax+vd(a)$, which is impossible. Thus $vx\not \in x\Tn$ and
$\Tn$ is not right duo, provided $R$ is commutative.
\end{proof}

\begin{remark}\label{commutative remark}
 One can easily check, with the help of Theorem \ref{th description of
 derivations}(2)(b), that if $R$ and $\de\ne 0$ are as in the proof of
 the implication $(a)\Rightarrow (b)$ in the above theorem, then
 $\Tn$ satisfies all necessary conditions from
 Proposition \ref{necessary conditions} for an Ore extension to
 be right duo.
\end{remark}
\begin{example}
Let $K=L(x_i\mid i=0,1,\ldots)$, $A=K[x]$, $\phi$ and $R=A\oplus
M$ be as in Example \ref{example 1}. Then $\phi (A)\subseteq
L(x_i\mid i=1,2,\dots)=\hat{L}$. Thus $x_0\in K$ is transcendental
over $\hat{L}$ and Theorem \ref{theorem duo rings}(2) implies that
$R[x;\si,\de_{x_0}]$ is a right duo ring.
\end{example}

It is easy to construct an example of a field $K$ with an
endomorphism $\phi$ such that all possibilities from Corollary
\ref{when ore ext is duo} occur.
\begin{example}Let $R=K\oplus vK$, where
  $K=F(X)$ is a field of rational functions over a field
 $F$ in the set $X=\{x_i\mid i=0,1,2,\ldots\}$ of indeterminates
 and $\phi\colon K\rightarrow K$ is an
 an $F$-linear homomorphism defined by
 setting $\phi(x_i)=(x_{i+1})^{i+1}$ for $x_i\in X$.

 It is
 easy to see that $x_0$ is transcendental over $\phi(K)$ as
 $\phi(K)\subseteq F(X\setminus\{x_0\})$, while $x_i$ is
 algebraic over $\phi(K)$ of degree $i$ for all $x_i\in
 X\setminus\{x_0\}$.

 Thus, by  Corollaries \ref{when ore ext is duo}, \ref{case A is a
 field} and Lemma \ref{outer derivations}, respectively,
 we have:
\begin{enumerate}
  \item $R[x;\si,\de_{x_0}]$ is a right duo ring.
  \item If $k\geq 1$, then $(x^k-(x_k)^k)R[x;\si,\de_{x_k}]$ is not a two-sided ideal of
 $R[x;\si,\de_{x_k}]$ but $fR[x;\si,\de_{x_k}]$ is a two-sided ideal, for all
 polynomials   $f\in R[x;\si,\de_{x_k}]$ with $\deg(f)<k$.
 \item  If $k\geq 2$, then $\de_{x_k}$ is an outer  $\si$-derivation
 of $R$.
 \end{enumerate}
\end{example}

One can  check  that, for any $k\geq 2$, the Ore extension
$T=R[x;\si,\de_{x_k}]$ from the above example satisfies all
necessary conditions from Proposition \ref{necessary conditions}
for the ring $T$ to be right duo. This means that conditions
obtained by Marks in \cite{Marks} for an Ore extension to be right
duo are not sufficient.

 We close the paper  formulating the
following problems:

\begin{Problem}
Let $B$ be a ring. Find necessary and sufficient conditions, in
terms of properties of $B,\;\tau$ and $\de$, for an Ore extension
$B[x;\tau,\de]$ to be a right duo ring.
\end{Problem}

\begin{Problem}
Let $B$ be a split corner  subring of $T=B\oplus M$ with $M^2=0$.
Find necessary and sufficient conditions, in terms of properties
of $B$ and the $(B-B)$-bimodule $M$, for the ring $T$ to be right
duo.
\end{Problem}

Theorems \ref{construction of duo rings} and \ref{theorem duo
rings} provide examples of right duo rings of the form $B\oplus M$
as above. Those are of the form described in Proposition \ref{main
prop}.

If $B$ is a commutative noetherian ring  then, as the following
proposition shows, $B[x;\tau,\de]$ is never a one-sided duo ring,
except the case $B[x;\tau,\de]=B[x]$. Nevertheless, by Remark
\ref{commutative remark}, there exist such noncommutative Ore
extensions which satisfy all necessary conditions from Proposition
\ref{necessary conditions}. By Proposition \ref{quasi-duo}, these
Ore extensions are quasi-duo rings.
\begin{prop}
 Let $B$ be a commutative noetherian ring.
 If the Ore extension $B[x;\tau,\de]$ is a right (left) duo ring,
 then $\tau=\mbox{\rm id}_B$ and $\de=0$, i.e. $B[x;\tau,\de]=B[x]$ is
 a commutative polynomial ring.
\end{prop}
\begin{proof}
 If $B[x;\tau,\de]$ is left duo, then the thesis is a consequence
 of Theorem \ref{Marks characterization of left duo}.

 Suppose that $B[x;\tau,\de]$ is a right duo ring which is
 noncommutative. Then, by Proposition \ref{nice corner}, there
 exists a noncommutative Ore extension $\Tn$ which is right duo,
 such that $R=A\oplus M$ is a split corner extension of a ring $A$
 with $M^2=0$, where $M=\ker\si\ne 0$ and $\si|_A=\mbox{id}_A$. Since
 $R$ is a factor ring of $B$, $R$ is commutative and noetherian.

 By
 Proposition \ref{reduction}, there is an ideal $J$ of $R$ such
 that $R/J\simeq A'\oplus M'$, where $M'$ is an
$(A'-A')$-bimodule which is simple as a right $A'$-module and
faithful as a left $A'$-module. Proposition \ref{necessary
conditions}(3) guarantees that $J$ is a $(\si-\de)$-stable, so
$\Tn/(J\Tn)\simeq (R/J)[x;\si,\de]$, where $\si$ and $\de$ denote
also the maps induced on $R/J$, i.e. replacing $R$ by $R'$,  we
may assume that the commutative ring $R=A\oplus M$, where $M$ is
simple as a right $A$-module and faithful as a left $A$-module.
Thus $R=A\oplus M$ is a ring considered in Sections 3 and  4. Now,
since $R$ is commutative,    Theorem \ref{theorem duo rings}(2)
yields that $\Tn$ is not right duo. This contradicts our
assumption and completes the proof of the proposition.
\end{proof}

\bf Acknowledgments. \rm I thank Greg Marks for careful reading
the manuscript and for simplifying the arguments used in the proof
of Proposition 4.5.
\end{document}